\documentclass[aps,reprint]{revtex4}
\usepackage{epsfig}
\usepackage{amsmath}
\usepackage{epsfig}

\begin{document}
\title
{A solution to Ahmed's Integral(II)}    
\author{Dona Ghosh,~M.Sc.~(Maths, IIT-Kharagpur India)\\
B3-11, 304, Astavinayak Society, Sector-4, Vashi, Navi-Mumbai 400703, India\\
rimidonaghosh@gmail.com}
\begin{abstract}
In the year 2000, Ahmed proposed a family of  integrals in the American Mathematical Monthly which invoked a considerable response then. Here I would like to present  another solution to this family of integrals. I propose to call this as Ahmed's Integral (II) in the light of the
well known Ahmed's Integral.
\end{abstract}
\maketitle
Earlier, an integral proposed by Ahmed in 2001-2002 [1,5] has been well discussed in the books [2-4], and included in mathematical encyclopedias and dictionaries. A Google search by ``Ahmed's Integral'' brings over 70 hits to view.

Fascinated by the popularity  of Ahmed's Integral, I found that Ahmed has proposed one more interesting integral in the year 2000 [6,7]. This integral also received a good response, when 31 authors and two problem solving groups proposed its solutions. The solution of Peter M. Jarvis (Georgia) was published entitled `A family of Integrals'[7]. 

The family of integrals proposed by Ahmed [6,7] is given as
\begin{equation}
I_{m,n}=\int_{0}^{\infty} \frac{d^m}{dx^m} \left(\frac{1}{1+x^2}\right)  \frac{d^n}{dx^n} \left(\frac{1}{1+x^2}\right) ~ dx. 
\end{equation}
In the following I wish to present another solution to (1). Let us define $f_n(x)=\frac{d^n}{dx^n} \left(\frac{1}{1+x^2}\right),n=0,1,2,3,...$, $f_0(x)=f(x)$ and $f(x)=
\frac{1}{1+x^2}$ and rewrite Eq.(1) as
\begin{equation}
I_{m,n}=\int_{0}^{\infty} f_m(x) f_n(x) dx,
\end{equation}
where $m$ and $n$ are non-negative integers.
We find that $f(x)$ satisfies a second order ordinary linear differential equation.
\begin{equation}
(1+x^2)f_2(x)+4xf_1(x)+2f(x)=0.
\end{equation}
Using Lebnitz rule, $n$ times differentiation of Eq. (3) w.r.t. $x$ gives us 
\begin{equation}
(1+x^2)f_{n+2}+2(n+2)xf_{n+1}+(n+2)(n+1)f_n=0.
\end{equation}
So we find that $f_n(\infty)=0$ and $f_n(0)=-n(n-1)f_{n-2}(0)$. One can check that $f_0(\infty)=0, f_1(\infty)=0, f_0(0)=1$ and $f_1(0)=0$. Hence, we get
\begin{equation}
f_n(\infty)=0,\quad f_n(0)= n! \cos(n\pi/2).
\end{equation}
Without a loss of generality, we can assume $m>n$.\\
{\bf Case 1:} $m+n=$ even\\
Integration by parts of Eq. (1) by treating $f_m(x)$ as first and $f_n(x)$ as second function  yields
\begin{equation}
I_{m,n}= \int_{0}^{\infty} f_m (x) f_n(x) dx= -f_{m}f_{n-1}-I_{m+1,n-1},
\end{equation}
Here onwards the argument of $f_i$ is $0$. The repeated use of the recurrence relation gives 
\begin{equation}
I_{m,n}= \epsilon_{m,n} I_{m+n,0} =\epsilon_{m,n} \int_{0}^{\infty} f_{m+n}(x) f_0(x) dx.
\end{equation}
Here $\epsilon_{m,n}=1$, when $m=n=$ even or when $m\ne n$. $\epsilon_{m,n}=-1$, when $m=n=$ odd. 
We now use a representation of $(1+x^2)^{-1}$ as
\begin{equation}
\frac{1}{1+x^2} = \int_{0}^{\infty} e^{-z} \cos xz~ dz.
\end{equation}
to express $f_{m+n}(x)$ and write
\begin{equation}
I_{m,n}=\epsilon_{m,n} \int_{0}^{\infty} \left (\frac{1}{1+x^2}
\int_{0}^{\infty}  z^{m+n} e^{-z} \cos xz ~dz\right) dx
\end{equation}
Next by using the cosine-Fourier transform of $(1+x^2)^{-1}$ i.e.,
\begin{equation}
\int_0^{\infty} \frac{\cos xz}{1+x^2} dx={\pi \over 2} e^{-z},
\end{equation}
we find 
\begin{equation}
I_{m,n}=\frac{1}{2}\pi \epsilon_{m,n}  \int_{0}^{\infty} z^{m+n} e^{-2z} dz= \epsilon_{m,n} \frac{(m+n)! \pi} {2^{m+n+2}},
\end{equation}
when $m+n$ is even.\\
{\bf Case 2:} When $m+n=$ odd\\
Again we assume $m>n$ without a loss of generality. The successive integration of $I_{m,n}$ by parts $(m-n)$ times, treating $f_{n}(x)$ as first and $f_{m}(x)$ as second function leads to
\begin{equation}
I_{m,n}=(-1)^m\sum_{j=1}^{m-n} f_{m-j} f_{n+j-1}-I_{m-(m-n),~n+(m-n)}.
\end{equation}
Next noting that $I_{m,n}=I_{n,m}$, we find 
\begin{eqnarray}
%\begin{equation}
I_{m,n}=\frac{(-1)^m}{4} \sum_{j=1}^{m-n} \left( \sin[(m+n) \frac{\pi}{2}]-(-1)^j\sin[(m-n) \frac{\pi}{2}] \right) (m-j)! (n+j-1)!~~
%\end{equation}
\end{eqnarray}
Two interesting cases arise here. These are
\begin{equation}
I_{2k+2,2k+1}=0 \quad \mbox{and} \quad I_{2k+1,2k}=-\frac{1}{2} [(2k)!]^2
\end{equation}
{\bf References:}\\
\noindent
$[1]$ Z. Ahmed, Amer. Math. Monthly, Problem Proposal 10884,
{\bf 108}(2001) 566; {\bf 109}(2002)  670-671. \\
$[2]$ J. M. Borwein, D.H. Bailey, and R. Girgensohn,  `Experimentation in Mathematics: Computational Paths to Discovery' (Wellesley, MA: A K Peters) (2004) pp. 17-20.\\ 
$[3]$ P. J. Nahin, `Inside Interesting Integrals' (Springer: New York) (2014) pp.190-194.\\
$[4]$ G. Boros and V.H. Moll, `Irresistible Integrals'
(Cambridge University Press: UK, USA) (2004) p.  277. \\
$[5]$ Z. Ahmed, `Ahmed's Integral: the maiden solution',
arxiv:1411.5169 [Math-HO]\\
$[6]$ Z. Ahmed, Amer. Math. Monthly, Problem Proposal 10777 {\bf 107} (2000) p. 83.\\
$[7]$ Z. Ahmed, `A family of Integrals', Amer. Math. Monthly {\bf 107} (2000) 956-957.
\end{document}